\input amstex
\documentstyle{amsppt}
\document
\magnification=1200
\NoBlackBoxes
\nologo
\pageheight{18cm}
\voffset =1.5cm
\def\<{\langle}
\def\>{\rangle}
\def\op{\operatorname}


\bigskip

\centerline{\bf MODULES AND MORITA THEOREM FOR OPERADS}

\medskip

\centerline{\bf M.~Kapranov${}^1$, Yu.~Manin${}^2$}

\medskip

\centerline{\it ${}^1$ University of Toronto,  Canada}

\smallskip

\centerline{\it ${}^2$Max--Planck--Institut f\"ur Mathematik, Bonn, Germany}

\bigskip

\centerline{\it Dedicated to A.~N.~Tyurin on the occasion of his 60th birthday}

\medskip

\centerline{\bf \S 0. Introduction and summary.} 

\medskip

{\bf (0.1) Morita theory.}  Let $A, B$ be two commutative rings.
If their respective categories of modules are equivalent,
then $A$ and $B$ are isomorphic. This is not  true anymore
if $A$ and/or $B$ are not assumed to be commutative.
Morita theory describes the relashionship between
$A$ and $B$ in this case. The fact that a ring is not determined
uniquely by its category of modules has deep implications
for non--commutative geometry which tends to substitute
an elusive non--commutative space by the category of sheaves
on it.

\smallskip

In this note we construct a fragment of Morita theory for operads.
More precisely, let $k$ be a field. For a $k$--linear operad $\Cal{P},$
we construct matrix operads $\roman { Mat}\,(n,\Cal{P})$
and prove that their respective categories of representations
(that is, algebras) are equivalent. 
\smallskip

In order to compare the situation with the classical one,
let us remind the exact statement of the Morita theorem.
For $A, B$ as above, denote by $A$--$\roman { Mod}$,
$B$--$\roman { Mod}$ their respective categories of 
left modules. They are equivalent iff one can find
a $(B,A)$--bimodule $M$ which is finitely generated
and projective as $B$--module. Then the functor
$A$--$\roman { Mod}$ $\to$ $B$--$\roman { Mod}:$
$$
N\mapsto M\otimes_A N
\leqno(0.1.1)
$$
establishes an equivalence,
and $A$ is isomorphic to an algebra of the type $e\,\roman { Mat}(n,B)e$
where $e$ is the idempotent defining $M$.

\smallskip

A $k$--linear operad $\Cal{P}$ can be considered as
an ``associative ring'' (or rather monoid) in the monoidal category 
$\bold{S}$--$\roman{Vect}$
whose objects are families of representations of all the
symmetric groups $\bold{S}_n, n\ge 0.$
The plethysm monoidal product in this category (denoted $\circ$ below) is not
symmetric. This adds a new dimension  of non--commutativity to
the situation. In particular, the notions of
left and right modules become asymmetric, and
in our Morita theorem we replace left modules
by $\Cal{P}$--algebras, whereas right modules remain
right $\Cal{P}$--modules in $\bold{S}$--$\roman{Vect}$. 
Denote the categories of these objects
$\Cal{P}$--$\roman{Alg}$ and $\roman{Mod}$--$\Cal{P}$ respectively.

\smallskip

Our argument proceeds as follows.

\medskip

(i) We use  a kind of relative tensor product operation
$\circ_{\Cal P}: \roman { Mod}$--$\Cal{P}
 \times \Cal{P}$--$\roman { Alg}\to \roman { Vect}$
(see sec. 1.3 below) to construct an operadic version of the functor 
(0.1.1).
  Let $\Cal{Q}$ be another $k$--linear operad.
Consider an $\bold{N}$--graded $\Cal{Q}$--algebra
$M$ which is simultaneously a right $\Cal{P}$--module
such that both structures are compatible in the sense
that will be made explicit in sec. 1.4. Then for any $A$
in $\Cal{P}$--$\roman{Alg}$ the product $M\circ_{\Cal{P}}A$
is in $\Cal{Q}$--$\roman{Alg}$.  
\smallskip

(ii) Given an $M$ in $\roman{Mod}$--$\Cal{P}$, we  can construct
its endomorphism operad $\Cal{Q}=\roman { Op\,End}_{\Cal{P}}(M)$. It consists
of the part of $\bigoplus_n \roman { Hom}\,(M^{\otimes n},M)$ compatible
with the grading and the action of $\Cal{P}.$
Then $M$ becomes a $\Cal{Q}$--algebra, so that the 
construction of (i) provides the functor
$\Cal{P}$--$\roman{Alg}\to \Cal{Q}$--$\roman{Alg}$.

\smallskip

(iii) Finally, if we take for $M$ a free module of rank $n$
in $\roman{Mod}$--$\Cal{P}$, its endomorphism operad
is denoted $\roman { Mat}\,(n,\Cal{P})$, and the functor produced
in (ii) turns out to be an equivalence of categories: see (1.8.1).

\medskip 
 
Hopefully, this construction can be extended to
a fuller statement providing the necessary conditions
for equivalence as well.

\smallskip

Moreover, the multifaceted analogies between
linear operads and rings suggest several
other Morita--like contexts. For example,
one can ask which equivalences may exist between the categories
of right modules for different operads.
Notice that equivalences between the categories of modules over algebras
of various operadic types can be studied in principle by
passing to the universal enveloping algebras and
applying Morita theory for associative algebras.

\medskip

{\bf (0.2) Plan of the paper.} The first section is devoted
to the Morita theorem. We first remind the generalities
on operads, operadic modules, and algebras, in order to
fix notation. We proceed by describing several examples
and constructions related to the operadic modules
which deserve to be better known: see e.~g. Theorem (1.6.3)
for the construction of a Lie algebra, which provides
in particular a family of canonical vector fields
on any formal Frobenius manifold. Finally, we prove
the operadic Morita theorem described above.

\smallskip

The simplest application of the classical Morita
theory to supergeometry furnishes the following fact: {\it the
category of $D$--modules on a supermanifold $M$
is equivalent to the one on its underlying
manifold $M_{\roman{red}}$}. (To get a more sensible
statement, $\bold{Z}_2$--grading must be introduced; for more
sophisticated generalizations see [P]). In fact,
the algebra of differential operators
of the $\bold{Z}_2$--graded exterior algebra $W:=\wedge_k(V)$ of a finite 
dimensional
vector space $V$ is isomorphic to the matrix
algebra $\roman{End}_k(W).$ Therefore, sheaves
of algebras $\roman{Diff}_M$ and $\roman{Diff}_{M_{\roman{red}}}$ 
are Morita equivalent.
In the second section of this paper we discuss
some operadic versions of this remark.

\medskip

{\it Acknowledgements.} The initial draft of this note was written 
in June 1995
when the first author was visiting
  Max--Planck Institut f\"ur
Mathematik.  He would like to thank the MPI
for hospitality and financial support. Since 1995,
some of the constructions of the first version were also
discovered by other people. In preparing the present version,
we tried to give ample references. The first author
would like to thank E.~Getzler for poining out the recent
work on pre-Lie algebras [CL] [Ku]. 
The second author would like to thank A.~Davydov for illuminating
correspondence about the plethysm product.
We are grateful to the referee for  several
very useful remarks on the earlier version of the paper. 
 
\newpage

\bigskip

\centerline{\bf \S 1. Modules over operads}

\smallskip

\centerline{\bf and abstract Morita theory.}

\bigskip

{\bf (1.1) Monoidal structures on S--{Vect}.}
We work over a fixed base field $k$ and denote by Vect the
category of $k$--vector spaces.
We denote by {\bf S} the category whose objects are the standard
finite sets $\<n\>=\{1, 2, ..., n\}$, $n\geq 0$ and morphisms are
bijections, i.e., $\text{Hom}_{\bold S}(\<n\>, \<n\>)=\bold{S}_n$
is the symmetric group. By an {\bf S}--space we mean
a contravariant functor ${\bold S}\to\text{Vect}$, i.e.,
a collection $\Cal V=\{\Cal V(n)\}_{n\geq 0}$ where
$\Cal V(n)$ is a vector space with a right $\bold{S}_n$-action.
As well known, an {\bf S}-spave $\Cal V$ defines a functor on the category of
all finite sets and their bijections. More precisely, if $I$ is
any set with $n$ elements, then we put
$$\Cal V(I) = \biggl( \bigoplus_{\phi: \<n\> \longrightarrow
I} \Cal V(n)\biggr)_{{\bold S}_n}.\leqno (1.1.1)$$ 
Here $\phi$ runs over all bijections.
Let {\bf S}--Vect be the category of {\bf S}--spaces.
It possesses three important monoidal structures.
Let $\Cal V=\{\Cal V(n)\}$ and $\Cal W=\{\Cal W(n)\}$ be two
{\bf S}--spaces. Set
$$(\Cal V\otimes\Cal W)(n) = \Cal V(n)\otimes_k\Cal W(n),\leqno
(1.1.2)$$
$$(\Cal V\boxtimes\Cal W)(n) = \bigoplus_{i+j=n} \text{Ind}_{\bold{S}_i
\times \bold{S}_j}
^{\bold{S}_n}(\Cal V(i)\otimes_k \Cal W(j)) = \bigoplus_{I\sqcup J = \<n\>}
\Cal V(I)\otimes_k \Cal W(J), \leqno (1.1.3)$$
$$\Cal V\circ\Cal W = \bigoplus_{n\geq 0} \Cal V(n)\otimes_{\bold{S}_n} \Cal W
^{\boxtimes n}.\leqno (1.1.4)$$
In other words, 
$$(\Cal V\circ\Cal W)(n) = \bigoplus_{p\geq 0}
\Cal V (p)\otimes_{\bold S_p}
\left(\bigoplus_{I_1\sqcup ... \sqcup I_p = \<n\>}  \Cal W(I_1)\otimes_k ...
\otimes_k \Cal W(I_p)
\right)$$
Here the ${\bold S}_p$-action on the right takes the summand corresponding to
$(I_1, ..., I_p)$ to the summand corresponding to $(I_{s(1)}, ..., I_{s(p)})$,
$s\in {\bold S}_p$. This summand will be different unless $I_\nu = I_{s(\nu)} = \emptyset$ for all $\nu$ not fixed by $s$. This means that

$$  (\Cal V\circ\Cal W)(n) = 
\bigoplus_{p\geq 0}\bigoplus_{a_1+...+a_p=n}
\biggl(\Cal V (p)\otimes_k \Cal W(a_1)\otimes_k ... \otimes_k \Cal W(a_p)\biggr)_{\op{Aut}(O(a_1, ..., a_p))}.
\leqno (1.1.5)$$
Here $O(a_1, ..., a_p)\subset \<p\>$ is the set of $\nu$ such that $a_\nu=0$
and $\op{Aut}(O(a_1, ..., a_p))\i {\bold S}_p$ is the group of self-bijections of this set.

Also, for $\Cal V\in {\bold S}$--Vect and $X\in$ Vect 
define the vector space
$$\Cal V\langle X\rangle = \bigoplus_{n\geq 0} \Cal V(n)
\otimes_{\bold{S}_n} X^{\otimes n}.\leqno (1.1.6)$$
The products (1.1.3--4) are analogous to the operations of
multiplication and composition of formal power series 
$\sum v_n x^n/n!$, and (1.1.6) is analogous to the evaluation
of such series.
The product (1.1.4) is known as plethysm. 
Notice that for a  ``constant" {\bf S}-space,
i.e., an {\bf S}-space  $\Cal W$ of the form
$\Cal W(0)=X$, $\Cal W(n)=0, n>0$ ,
the plethysm $\Cal V\circ \Cal W$ is again ``constant",
corresponding to the vector space $\Cal V \<X\>$.

Each of the three products makes {\bf S}--Vect into a $k$--linear
monoidal category, with $\otimes$ and $\boxtimes$ giving in
fact symmetric monoidal structures, whereas $\circ$ is non--symmetric.
The unit objects with respect to the three structures are:
$$\Cal Com : \quad \Cal Com(n)=k, n\geq 0, \leqno (1.1.7)$$
$${\bold 1} : \quad {\bold 1}(0)=k, \, {\bold 1}(n)=0, n\neq 0,
\leqno (1.1.8)$$
$${\bold I} : \quad {\bold I}(1)=k, \, {\bold I}(n)=0, n\neq 1.
\leqno (1.1.9)$$
Note that we have canonical identifications
$$(\Cal V\circ \Cal W)\langle X\rangle =
\Cal V\langle \Cal W\langle X\rangle\rangle,
\leqno (1.1.10)$$
$$(\Cal V\boxtimes\Cal W)\circ \Cal X = (\Cal V\circ \Cal X)\boxtimes
(\Cal W\circ \Cal X), \leqno (1.1.11)$$
$$(\Cal V\boxtimes \Cal W)\langle X\rangle =
\Cal V\langle X\rangle \otimes_k \Cal W\langle X\rangle,
\leqno (1.1.12)$$
which correspond to the familiar rules of dealing with power series. 
In addition, we have the natural morphisms of {\bf S}--spaces
$$(\Cal V_1\otimes \Cal V_2)\circ (\Cal W_1\otimes \Cal W_2)
\to (\Cal V_1\circ \Cal W_1)\otimes (\Cal V_2\circ \Cal W_2)
\leqno (1.1.13)$$
and of vector spaces
$$(\Cal V_1\otimes\Cal V_2)\langle X_1\otimes X_2\rangle\to
\Cal V_1\langle X_1\rangle \otimes \Cal V_2\langle X_2\rangle.
\leqno (1.1.14)$$
For example, (1.1.13) is the ``diagonal map"
sending the summand 
$$\biggl(\Cal V_1(p)\otimes\Cal V_2(p) \otimes\Cal W_1(a_1)\otimes\Cal W_2(a_1)
\otimes... \otimes \Cal W_1(a_p)\otimes \Cal W_2(a_p)\biggr)_{
\op{Aut}(O(a_1, ..., a_p))}
\subset$$
$$\subset \bigl((\Cal V_1\otimes \Cal V_2)\circ(\Cal W_1\otimes\Cal W_2)
\bigr)(n), \quad a_1+...+a_p=n,$$
see (1.1.5), into the tensor product of the summand
$$\biggl(\Cal V_1(p)\otimes\Cal W_1(a_1)\otimes ... \otimes\Cal W_1(a_p)\biggr)_
{\op{Aut}(O(a_1, ..., a_p))}
\subset (\Cal V_1\circ\Cal W_1)(n)$$
and of the summand
$$\biggl(\Cal V_2(p)\otimes\Cal W_2(a_1)\otimes ... \otimes\Cal W_2(a_p)\biggr)_
{\op{Aut}(O(a_1, ..., a_p))}
\subset (\Cal V_2\circ\Cal W_2)(n).$$
This map is injective, if the {\bf S}-spaces $\Cal W_1$ and $\Cal W_2$
are null in degree 0. The morphism (1.1.14) is defined
in a similar way: it is the degree 0 component of (1.1.13)
for the ``constant" {\bf S}-spaces $\Cal W_i$ corresponding to $X_i$. 

\bigskip

{\bf (1.2) Operads and modules.} 
A $k$-linear operad $\Cal{P}$ can be defined 
as a monoid object in 
($\bold{S}$--$\roman { Vect},\,\circ$). This means that $\Cal{P}$
is an $\bold{S}$--space endowed with associative multiplication
(which is a morphism in $\bold{S}$--Vect)
and a distinguished element
$$
\mu_{\Cal P} :\, \Cal{P}\circ\Cal{P}\to\Cal{P},\  1\in\Cal{P}(1)
\leqno(1.2.1)
$$
with the usual properties. Spelling these data out
with the help of (1.1.5), we get the maps
$$
\mu_{m_1,\dots ,m_l}:\,
\Cal{P}(l)\otimes\Cal{P}(m_1)\otimes\dots\otimes\Cal{P}(m_l)
\to\Cal{P}(m_1+\dots +m_l)
\leqno(1.2.2)
$$
satisfying May's axioms (see [M], [S], [Gi--K]).

\smallskip

A right, resp. left $\Cal P$--module is an {\bf S}--space $M$ together with
a multiplication morphism
$$
{}^M \mu: M\circ \Cal P\to M, \quad \text{resp.}\quad
\mu^M: \Cal P\circ M\to M
\leqno (1.2.3)
$$
satisfying the usual asociativity and unit requirements (cf.
[Mar1-2], [Re],  [BJT]). We denote by Mod--$\Cal P$, resp. $\Cal P$--Mod, the 
categories
of right, resp. left $\Cal P$-modules. Note that Mod--$\Cal P$
is an abelian $k$-linear
category, because the plethysm product $\circ$ is linear in the
first argument, but $\Cal P$--Mod is not even additive. 

\smallskip

Similarly, a $\Cal P$--algebra is a vector space $A$ together with
a morphism of vector spaces
$$\mu^A: \Cal P\langle A\rangle \to A\leqno (1.2.4)$$
satisfying the usual requirements. We denote by $\Cal P$--Alg the
category of 
$\Cal P$--algebras. It is also not additive. 

\smallskip

If $V$ is a vector space, then $\Cal P\langle V\rangle$ is a
$\Cal P$--algebra called  {\it the free algebra generated by $V$}.

Given an operad $\Cal P$, one defines (see [Ad], \S 2.3)
 the corresponding category
of operators (or PROP) $U(\Cal P)$ as follows. Objects of $U(\Cal P)$
are the same as for {\bf S}, i.e., the sets $\<n\>$, $n\geq 0$, while
$$\op{Hom}_{U(\Cal P)} (\<m\>, \<n\>) = \bigoplus_{f: \<m\>\to\<n\>}
\bigotimes_{i=1}^n \Cal P(f^{-1}(i)),\leqno (1.2.5)$$
where $f$ runs over all maps of sets $\<m\>\to \<n\>$ and
the value of $\Cal P$ on a set is defined by (1.1.1). The composition
is induced by the composition in $\Cal P$. The following
is then obvious by construction.

\proclaim {(1.2.6) Proposition}
The category of right $\Cal P$-modules is equivalent to that of
 contravariant functors
$U(\Cal P)\to\op{Vect}$.
\endproclaim

\medskip

\noindent {\bf (1.3) Examples of operads and modules.}

\vskip .1cm

\noindent {\bf (1.3.1) The endomorphism operad.} 
Any vector space $V$ determines the operad $\roman { Op\,End}\,(V)$
 with
$$
\roman { Op\,End}\,(V)(n) = \roman { Hom}_k\,(V^{\otimes n},V)
\leqno 
$$
and the multiplication law (1.2.2) given by the composition
of multilinear maps. A structure of the $\Cal P$--algebra on $V$
is the same as a morphism of operads $\Cal P\to \text{Op\,End}(V)$. 

 An example of a right $\text{Op\,End}(V)$--module is given by 
$V^\star=\{(V^*)^{\otimes n}\}_{n\geq 0}$. The module structure
is given by taking the superposition of an $n$--linear form $V^{\otimes n}\to 
k$
with  $n$ multilinear maps $V^{\otimes a_i}\to V$, $i=1, ..., n$. 

\smallskip

More  generally, if $(\Cal C, \otimes)$ is any
$k$-linear
symmetric monoidal category, any object $V\in\Cal C$ gives rise
to an operad $\text{Op\,End}_{\Cal C}(V)$ defined
in a similar way. 

\smallskip

\noindent {\bf (1.3.2) The trivial operad.} 
The unit object ${\bold I}\in ({\bold S}$--Vect, $\circ)$
is an operad called the trivial operad. An {\bf I}--algebra is the same
as a vector space and an {\bf I}-module (left or right)
is the same as an {\bf S}--space. 

\smallskip

\noindent {\bf (1.3.3) Free and projective right modules.} 
 Any operad $\Cal P$ is a left and a right module over itself. 
It follows that for any vector space $V$ the {\bf S}--space
$V\otimes\Cal P=\{V\otimes_k \Cal P(n)\}_{n\geq 0}$ is a right
$\Cal P$--module. Such modules will be called free (of finite
rank, if $\dim\,V<\infty$).  A right module will be called
projective, if it is a direct summand of a free module.

For example $\text{Op\,End}(V)$ is isomorphic,
as a right module over itself,  to the direct sum of $\dim(V)$ copies
of $V^\star$, so $V^\star$ is a projective module.

\smallskip

\noindent {\bf (1.3.4) The commutative operad.}  
 The {\bf S}--space $\Cal Com$ has a natural operad structure with
all the maps (1.2.2) being the canonical identifications
$k\otimes ...\otimes k\to k$. A $\Cal Com$--algebra is the same
as a commutative algebra with unit in the usual sense. Note that
$\Cal  Com = \text{Op\,End}(k)$. 

  Let $F$ be the category whose objects are the  finite
sets $\<n\>$, $n\geq 0$ as above, but morphisms are all maps
of finite sets. Thus {\bf S} is the subcategory of
$F$ formed by all objects and all isomorphisms. By an $F$--space,
 we mean a contravariant functor $\Cal V: F\to \text Vect$,
whose value on $\<n\>$ will be denoted $\Cal V(n)$.
Such a functor gives, in particular, an {\bf S}--space.
We denote by $F$--Vect the category of $F$--spaces.

The category $F$  should not be confused
with the category $\Gamma$ of Segal [Se]
whose objects are finite {\it pointed} sets $[n] = \{0, 1, ..., n\}$,
$n\geq 0$ and morphisms are all maps taking 0 to 0.

\proclaim {(1.3.5) Proposition} The category $\roman{Mod}$--$\Cal Com$
of right $\Cal Com$--modules is equivalent to $F$--$\roman{Vect}$
so that the forgetful functor $\roman{Mod}$--$\Cal Com\to
{\bold S}$--$\roman{Vect}$ is identified with the restriction functor
$F$--$\roman{Vect}\to {\bold S}$--$\roman{Vect}$.

\endproclaim

\noindent {\sl Proof:} Follows from the definition (1.2.5) of
the category $U(\Cal P)$ and Proposition 1.2.6. 

\smallskip

For example, if $A$ is a commutative algebra, then
the morphism of operads $\Cal Com\to \text{Op\, End}(A)$ together
with the construction of (e)
makes $A^\star$ into a right $\Cal Com$-module. On the other hand,  $\<n\>\mapsto
A^{\otimes n}$ is a covariant functor $F\to \roman{Vect}$
(cf. [L], Proposition 3.2 or [Pi], Sect. 1.7), and thus
$\<n\>\mapsto (A^*)^{\otimes n}$ is an $F$-space.

 \medskip

\noindent {\bf (1.3.5) The associative operad and simplicial spaces.}
 Let $\Cal Ass$ be the operad whose algebras are associative algebras with
unit ([Gi--K]). Thus $\Cal Ass(n)$ is the regular representation of 
$\bold{S}_n$.
Let us describe right $\Cal Ass$--modules. Denote by $\widetilde F$
the category whose objects are the finite sets $\,n\>$ as before
and a morphism $\Phi: \<m\>\to [n]$ consists, first, of a map
$\phi: \<m\>\to \<n\>$ in the ordinary sense and, second, of a choice
of a total order on each fiber $\phi^{-1}(i)$. The composition of
morphisms in $\widetilde F$ is defined using the lexicographic ordering
of the fibers of a composition. Thus we have a functor 
$\widetilde F\to F$. 
 As before, we define a 
$\widetilde F$--space  as a contravariant functor $\widetilde F\to
\text{Vect}$. The proof of the following  proposition is also straightforward
from (1.2.5-6).

\proclaim { (1.3.6) Proposition} The category of right $\Cal Ass$--modules
is equivalent to the category of $\widetilde F$--spaces. 
\endproclaim

Let us now discuss the relation between  $\widetilde F$-spaces and simplicial spaces.
Let $\Delta$ be  the standard simplicial category, with objects
$[n] = \{0, ..., n\}$ and monotone maps as morphisms.
Thus contravariant functors $V_\bullet: \Delta\to\text{Vect}$ are simplicial
vector spaces, with the value of $V_\bullet$ at $[n]$ denoted $V_n$. 
We will use the standard notation for the  face and degeneracy operators
$$\partial_{n,i}: V_n\to V_{n-1},\, 0\leq i\leq n, \quad
s_{n,i}: V_n\to V_{n+1}, \, 0\leq i\leq n.\leqno (1.3.7)$$
On the other hand, let $\Delta_+\subset F$ be the subcategory with the
same objects $\<n\>$ bu only monotone maps as morphisms. 
The bijections $[n]\to \<n+1\>$ (taking $i$ to $i+1$) identify
 $\Delta$ with the
full subcategory of $\Delta_+$ on objects $\<n\>, n>0$.
Thus a $\Delta_+$-space $\Cal V$ is the same as, first, a simplicial
space $V_\bullet$ with $V_n = \Cal V(n+1)$ and, second, the datum
of a vector space $V_{-1} = \Cal V(0)$ together with
a linear map $\partial_{0,0}: V_0\to V_{-1}$  satisfying
$\partial_{0,0} \circ\partial_{1,0} = \partial_{0,0}\circ\partial_{1,1}$. 
Such objects are traditionally called augmented simplicial
spaces, and $\partial_{0,0}$ is called the augmentation. 
Note that every simplicial space can be considered as an augmented one,
by taking $V_{-1}=0$.

\proclaim {(1.3.8) Proposition} (a) There is an embedding of categories
$\Delta_+\subset F$, identical on objects, so for an $\widetilde F$-space
$\Cal V$ the collection of $V_n = \Cal V(n+1)$, $n\geq -1$,
forms an augmented simplicial space. 
 
(b) If $V_\bullet =\{V_n\}_{n\geq -1}$ is an augmented simplicial space, then
the collection of $V_{n+1}\otimes_k k[\bold{S}_{n+1}]$ forms an $\widetilde F$--space. 
\endproclaim

For example, if $A$ is an associative algebra with 1, then $A^\star$ is a right
$\Cal Ass$-module. On the other hand, setting $V_n = (A^*)^{\otimes (n+1)}$
we get an augmented simplicial space with $\partial_{n,i}$ given by
contractions with 1 and $s_{n,i}$ by inserting the map $A^*\to A^*\otimes_k
A^*$ dual to the multiplication. 

\smallskip

\noindent {\sl Proof of (1.3.8):} (a) If $\phi: \<m\>\to \<n\>$ is a
monotone map, we take on each $\phi^{-1}(i)$ the total order
induced from $\<m\>$. This gives a morphism $\widetilde \phi$ of $\widetilde F$.
If $\phi$ and $\psi$ are composable monotone maps, then one sees readily
that $\widetilde{\phi\psi} = \widetilde \phi \widetilde\psi$. 

(b) This is a consequence of the following property of $\widetilde F$.
Let $\Phi = (\phi, \gamma)$ be any morphism of $\widetilde F$,
so $\phi: \<m\>\to \<n\>$ is a map of sets and $\gamma$ is a system
of total orders on the $\phi^{-1}(i)$. Then there is a unique permutation
$s\in\bold{S}_m$ such that $\Phi = \widetilde\phi s$ where
$\widetilde \phi$  is as in the proof of (a) and $s$
 is considered as a morphism of
$\widetilde F$ in an obvious way (its fibers are singletons
so do not need ordering). 

\medskip

\noindent {\bf (1.3.9) The stable curves operad and its modules.} 
For a finite set $I$ we denone by $\overline{M}_{g,I}$ the Deligne-Mumford
stack classifying stable
    curves  of genus $g$ with marked points $(x_i)_{i\in I}$ labeled by
$I$ cf. [Kn]. For any injective map $\phi: I\to J$ of finite sets
and ang $g$ such that $\overline{M}_{g, I}\neq\emptyset$ there is
a natural morphism of stacks $\overline{M}_{g, J}\to \overline{M}_{g, I}$
called the stable forgetting (cf. [Man2], p. 93).

\smallskip

We first consider the case $I=[n] = \{0,1, ..., n\}$. The point $x_0$ on a
curve from $\overline{M}_{g, [n]}$ will be called the root point. 
The group $\bold{S}_n$ acts upon $\overline{M}_{g, [n]}$ by renumbering
all points except for $x_0.$ Moreover,  we have the  morphisms
$$
\overline{M}_{g,[l]}\times \overline{M}_{0, [m_1]}\times \dots
\overline{M}_{0, [m_l]}\to \overline{M}_{g,[m_1+\dots +m_l]}
\leqno (1.3.10)
$$
gluing the  root point of the universal curve parametrized by
$\overline{M}_{0,[m_i]}$ to the $i$--th labeled point
of the universal curve parametrized by
$\overline{M}_{g,[l]}$, $i=1,\dots ,l.$

\smallskip

If $g=0$, we get an operad $\overline{M}$ in the monoidal category of smooth
projective manifolds with $\overline{M}(n) =  \overline{M}_{0,n+1}$,
 $n\ge 2$ and $\overline{M}(1) = \{\op{pt}\}$, $\overline{M}(0)=\emptyset$
The compositions (1.2.2) not involving $\overline{M}(1)$ are given by
  (1.3.10) while the unique element of $\overline{M}(1)$ is the unit.
To produce a $k$--linear operad we put $H_*\overline{M}(n):=
H_*(\overline{M}_{0,n+1},k)$ and define the structure
maps via pushforward of the geometric maps.  
\smallskip

Algebras over $H_*\overline{M}$ (more precisely, cyclic
algebras, cf. below)
are called Cohomological Field Theories.
Formal completion of such an algebra $H$ at zero has
a natural structure of a formal Frobenius manifold.
 The theory of Gromov--Witten
invariants produces such a structure on the cohomology
of any smooth projective algebraic manifold.
For all of this  see e.g. [Man2]. 

\medskip

We now construct a family of right modules over $\overline{M}$.
Let  $S$ be a finite set (possibly empty) and let $g$ be
such that $\overline{M}_{g, S}\neq\emptyset$. For any $n$ consider
the $\bold {S}_n$-action on $\overline{M}_{g, S\sqcup \<n\>}$
given by renumbering points with labels in $\<n\>$.
 We have then
the morphisms
$$
\overline{M}_{g, S\sqcup \<l\>}\times \overline{M}_{0, [m_1]}\times \dots
\overline{M}_{0, [m_l]}\to \overline{M}_{g, S\sqcup \<m_1+\dots +m_l\>},
\leqno (1.3.11)
$$
defined in a way completely analogous to (1.3.10). These morphisms
define on
 the {\bf S}-stack $\overline{\Cal M}^S_g$, $\overline{\Cal M}^S_g(n) = 
\overline{M}_{g, S\sqcup \<n\>}$,  the structure of right
$\overline{M}$-module. 
  Again, in order to pass
to the $k$--linear situation, it suffices to
apply any homology theory with coefficients in $k$.
One can even consider Chow groups because the K\"unneth
formula holds for the left hand side of (1.3.11)
in the Chow theory. 

\smallskip

Actually, both sides
of (1.3.11) admit compatible morphisms onto $\overline{M}_{g,S}:$
at the left hand side, project onto $\overline{M}_{g, S\sqcup \<l\>}$,
forget $x_1,\dots ,x_l$ and stabilize, at the right hand side
forget $x_1,\dots ,x_{{m_1}+\dots +{m_l}}$ and stabilize.
For a stable $S$-pointed curve $(C, (x_s))$ of genus $g$
let $\widehat{C}^{l}_{(x_s)}$ be the fiber of 
$\overline{M}_{g, S\sqcup \<l\>}\to \overline{M}_{g, S}$ over
the point represented by $(C, (x_s))$. This is in fact an algebraic
variety (not just a stack).
Now, the above discussion leads to the following result.

\proclaim {(1.3.12) Proposition} (a) For any stable $S$-pointed curve
$(C, (x_s))$ the collection of the $\bold{S}_l$-varieties
$\widehat{C}^{l}_{(x_s)}$ forms
a right module over the operad $\overline{M}$.

(b) If $\phi: T\to S$ is an injective map and $(D, (y_t))$
is the stable $T$-pointed curve obtained from $(C, (x_s))$
by stable forgetting, then we have  natural morphisms
of ${\bold S}_l$-varieties $\widehat{C}^{l}_{(x_s)}\to 
\widehat{D}^{l}_{(y_t)}$
which form a morphism of right $\overline{M}$-modules. 
\endproclaim

For example, if $S=\emptyset$ and $C$ is smooth, then $\widehat{C}^{l}$ is 
the  Beilinson--Ginzburg--Fulton--MacPherson
``resolution of diagonals'' of $C^l$, see  [BG]. 
The construction of [BG] is actually applicable to any smooth curve
(stable or not) and  the $\overline{M}$-module
structure in this case can also be constructed directly, using
Proposition 3.8 of {\it loc. cit.} One can also compare with [Mar2]
which essentially deals with a real version of $\widehat{C}^l$
(for $C$ a circle).

As before, one can produce from each geometric module in Proposition
1.3.12, a
$k$--linear homology module. Such modules for different choices
of $(C, (x_s))$
form a constructible sheaf
over   $\overline{M}_{g,S}$, smooth along the natural stratification
by the type of the dual graph.

\smallskip

The last remark concerns enlarged symmetry of $\overline{M}_{g,l+1}$.
In fact, the whole $\bold{S}_{l+1}$ acts upon this space rather than
its subgroup $\bold{S}_l.$ The axiomatization of this symmetry
leads to the notion of the {\it cyclic operad}
introduced and studied in [Ge--K]. Algebras and modules over
a cyclic operad may also admit a cyclic structure,
and we may ask for a cyclic version of Morita theory.
We leave this question for  future research.

\medskip

{\bf (1.4) Relative plethysm.} We will now review the
relative plethysm (or circle-over construction) for modules
over an operad ([Re]).

Let $\Cal P$ be an operad, $M$ a right $\Cal P$--module and $N$ a left
$\Cal P$--module. Their relative plethysm $M\circ_{\Cal P}N\in {\bold S}$--Vect
is defined as the cokernel of (the difference of) the two morphisms
$$\partial_0, \partial_1: M\circ \Cal P\circ N 
{\longrightarrow\atop\longrightarrow}  M\circ N, \quad \partial_0 = {}^M\mu
\circ \text{Id}_N, \,\, \partial_1 = \text{Id}_M\circ \mu^N.\leqno (1.4.1)$$
This construction is similar to the usual tensor product of
a right and a left module over an algebra and in fact can
be performed for ``modules" over
 any monoid object in any monoidal category
(provided the cokernels exist).
Similarly, let $A$ be a $\Cal P$--algebra. The relative evaluation
$M\circ_{\Cal P}\langle A\rangle \in$ Vect is, by definition, the cokernel of
the two morphisms
$$\partial_0, \partial_1: (M\circ \Cal P)\langle A\rangle 
{\longrightarrow \atop\longrightarrow} M\langle A\rangle,\quad
\partial_0 = {}^M\mu\langle\text{Id}_A\rangle, \,\,
\partial_1 =\text{Id}_M\langle\mu_A\rangle.\leqno (1.4.2)$$
The following canonical identifications are proved by mimicking the
standard arguments for modules over an algebra:
$$\Cal P\circ_{\Cal P}N =N, \quad M\circ_{\Cal P}\Cal P = \Cal P, \quad
\Cal P\circ_{\Cal P}\langle A\rangle = A.\leqno (1.4.3)$$

\medskip

{\bf (1.4.4) Example.}  Let $V$ be a finite--dimensional 
$k$--vector space and $\Cal P = \text{Op\, End}(V)$, $A=V$, $M=V^\star$,
see (1.3.1). Then we claim that 
$$V^\star\circ_{\text{Op\, End}(V)} \langle V\rangle = k.$$
Indeed, we have a morphism of vector spaces
$$\widetilde \phi:\quad  V^\star\langle V\rangle = \bigoplus_{n\geq 0} 
(V^*)^{\otimes n}
\otimes_{S_n} V^{\otimes n} \quad \to\quad k,
\leqno (1.4.5)$$
which on the $n$--th summand is induced by the $n$--th tensor power
of the canonical pairing. Clearly $\widetilde\phi$ descends to a surjective
morphism
$\phi: V^\star\circ_{\text{Op\,End}(V)} \langle V\rangle\to k$. Let us
prove that $\phi$ is injective. For this, notice first that the image
in $V^\star\circ_{\text{Op\,End}(V)} \langle V\rangle$ of the $n$th summand in
$V^\star\langle V\rangle$ factors through
$$\biggl((V^*)^{\otimes n} \otimes_{{\text End}(V)^{\otimes n}}V^{\otimes n}
\biggr)_{S_n} = k.$$
In other words, we have morphisms $\psi_n: k\to 
V^\star\circ_{\text{Op\,End}(V)} \langle V\rangle$
such that $\sum_{n\geq 0} \psi_n$ is surjective. Let us now prove that
the images of the $\psi_n$ coincide with each other.  For this,
choose an identification $V\to k^d$, $d=\dim\,V$ and use this to
make $V$ into a commutative associative algebra (the direct
sum of $d$ copies of $k$). We get an element $m\in\text{Op\, End}(V)(2)$
giving this algebra structure. Now, using the $\text{Op\, End}(V)$--linearity
conditions for $m$ and several copies of 1, we identify the
$\text{Im}(\psi_n)$ with each other. 

\bigskip

{\bf (1.5) Bimodules and functors between categories of operadic
algebras.} Let $\Cal{P},\,\Cal{Q}$ be two operads.
A $(\Cal{Q},\,\Cal{P})$--bimodule ([Mar1] [Re]) is, by definition, a
 space $M\in \bold{S}$--Vect  together with a  right $\Cal{P}$--module
structure  and 
 a left $\Cal{Q}$--module structure on ${M}$,
which commute, i.e., give rise to a well-defined map
$$
\Cal Q\circ M\circ\Cal P \to M. 
\leqno (1.5.1)
$$
 
\medskip 

{\bf (1.5.2) Examples.} (a) 
Any operad $\Cal{P}$ is a $(\Cal{P},\Cal{P})$-bimodule.

\smallskip

(b) Let $V,W$ be two vector spaces. Define the {\bf S}--space
$\text{Op\, Hom}(V,W)$ by
$$\text{Op\, Hom}(V,W)(n) = \text{Hom}_k(V^{\otimes n}, W).$$
Then $\text{Op\, Hom}(V,W)$ is an $(\text{Op\, End}(W), 
\text{Op\, End}(V))$--bimodule.

\smallskip

Similarly, if $(\Cal C,\otimes)$ is any symmetric monoidal category
and $V,W$ are objects  of $\Cal C$, we define an
 $(\text{Op\, End}_{\Cal C} (W), 
\text{Op\, End}_{\Cal C} (V))$--bimodule  $\text{Op\, Hom}_{\Cal C} (V,W)$

\smallskip

Let $L$ be a left $\Cal Q$--module, $N$ a right $\Cal P$--module and $M$,
as before, a $(\Cal Q, \Cal P)$--bimodule. Then $L\circ_{\Cal Q} M$ is
naturally a right $\Cal P$-module, $M\circ_{\Cal P}N$ is naturally
a left $\Cal Q$--module and we have a canonical associativity isomoprhism
$$(L\circ_{\Cal Q}M)\circ_{\Cal P}N \simeq L\circ_{\Cal Q}(M\circ_{\Cal P}N),
\leqno (1.5.3)$$
which allows us to write more complicated iterated
relative plethysm without parentheses. Similarly, if $A$ is a $\Cal 
P$--algebra,
we have a canonical isomoprhism
$$
(L\circ_{\Cal Q}M)\circ_{\Cal P}\langle A\rangle \simeq
L\circ_{\Cal Q}\langle M\circ_{\Cal P}\langle A\rangle\rangle.
$$
As a consequence, note a formula for the relative evaluation on a free
algebra:
$$
M\circ_{\Cal Q}\langle \Cal Q\langle A\rangle\rangle \simeq
M\langle A\rangle
\leqno (1.5.4)
$$
(take $\Cal P={\bold I}$ to be the trivial operad). 

\medskip

{\bf (1.5.5) Examples.} (a) Let $V$ be a finite--dimensional
vector space. Then $V\otimes\Cal Com = \{V\}_{n\geq 0}$ is naturally
an $(\text{Op\,End}(V), \Cal Com)$--bimodule. The arguments used in
Example 1.4.4 establish an identification
$$V^\star\circ_{\text{Op\, End}(V)}(V\otimes\Cal Com) \simeq \Cal Com.$$

(b) Let $V,W,X$ be finite--dimensional $k$--vector spaces. Then we have
an identification
$$\text{Op\, Hom}(W,X)\circ_{\text{Op\, End}(W)} \text{Op\, Hom}(V,W)
\simeq \text{Op\,Hom}(V,X)$$
(as bimodules). This follows easily from the above example and Example
1.4.4, because $\text{Op\, Hom}(W,X)$ is isomorphic, as a right $\text{Op\, 
End}
(W)$--module, to the direct sum of $\dim\,X$ copies of $W^\star$, and
$\text{Op\, End}(W)$ is isomorphic, as a right module over
itself, to the direct sum of $\dim\,W$ copies of $W^\star$. 

\medskip

The importance of $(\Cal Q, \Cal P)$--bimodules for us is that
any such bimodule $M$ defines a functor 
$f_M:\,\Cal{P}$--$\roman{Alg}\to
\Cal{Q}$--$\roman{Alg},$
$$  
f_M(A)=M\circ_{\Cal{P}}\langle A\rangle .
\leqno(1.5.6)
$$
The action of $\Cal{Q}$ on $f_M(A)$ is transferred from
the action on $M.$ Relative plethysm of bimodules corresponds to
the composition of functors: 
$$f_{N\circ_{\Cal Q}M} = f_N\circ f_M,\leqno (1.5.7)$$
as it follows from (1.5.4).

\medskip

{\bf (1.5.8) Example.} Let $\Cal P={\bold I}$ be the trivial operad,
so that $\Cal P$--$\text{Alg}=\text{Vect}$. An $({\bold I}, {\bold 
I})$--bimodule
$M$ is just an {\bf S}--space. The functor $f_M:\text{Vect}\to\text{Vect}$
is then
$$V\mapsto M\langle V\rangle = \bigoplus_{n\geq 0} M(n)\otimes_{
\bold{S}_n} V^{\otimes
n}.$$
Functors of this kind are called analytic in [J].

\smallskip

The following corollary of Example 1.5.5(b) is the first instance
of the general Morita--type theorem to be proved later.

\proclaim{\quad (1.5.9) Theorem}
Let $V, W$ be any two finite--dimensional $k$--vector spaces.
Then the functor $f_{\roman{Op\, Hom}\,(V,W)}$ induces an equivalence between
the categories of algebras over $\roman{Op\, End}\,(V)$ and over
$\roman{Op\, End}\,(W)$. In particular, all these categories are equivalent
to the category of commutative algebras, as $\Cal Com = \roman{Op\, 
End}\,(k)$.
\endproclaim

\smallskip

In fact, the functor $f_{\text{Op\, Hom}(V,W)}$ can be described in
a more explicit way, using the tensor product of modules over a ring.

\smallskip

\proclaim {\quad (1.5.10) Proposition} Let $V,W$ be as before and $A$ be an
$\roman{Op\, End}(V)$--algebra. Then, as a vector space,
$$
\roman{Op\,Hom}\,(V,W)\circ_{\roman {Op\, End}\,(V)}\langle A\rangle 
\simeq \roman{Hom}_k(V,W)\otimes_{\roman {End}\,(V)} A,
$$
where we view $A$ as a left module over the ring $\roman{End}(V) = 
\roman{Op\,End}(V)(1)$. 
\endproclaim

\smallskip

{\bf Proof.} As we pointed out already,
$\text{Op\, Hom}(V,W)$ is isomorphic, as a right $\text{Op\, End}(V)$--module,
to a direct sum of $\dim\,W$ copies of $V^\star$. Notice also that
 both sides of the proposed equality are additive in $W$. Thus,
it is enough to treat the particular case $W=k$, so that
$\text{Op\, Hom}(V,W)=V^\star$. Next, since the functor
$$
f_{\text{Op\, Hom}\,(k,V)}: \text{Op\, End}\,(k)\,\text{--Alg} = 
\Cal Com\,\text{--Alg} \to \text{Op\, End}\,(V)\,\text{--Alg}
$$
is an equivalence, we can assume that $A$ lies in the image of this
functor. However, $\text{Op\, Hom}(k,V)$ is isomorphic to $V\otimes_k\Cal Com$
as a right $\Cal Com$--module. Therefore, for a commutative algebra $B$
we have
$$f_{\text{Op\, Hom}(k,V)}(B) = (V\otimes_k\Cal Com)\circ_{\Cal Com}
\langle B\rangle =
V\otimes_k B.$$
Now, the left  and the right hand sides of the proposed equality are,
respectively
$$V^\star\circ_{\text{Op\, End}(V)} \langle V\otimes_k B\rangle, \quad
V^*\otimes_{\text{End}(V)} (V\otimes_k B).$$
But both these spaces are canonically identified with $B$: the first one
in virtue of Theorem 1.5.9 and Example 1.5.5, the second one  by elementary
linear algebra. Proposition is proved.  

\bigskip

{\bf (1.6) Right modules as a tensor category: the role of
the product $\boxtimes$.}
Consider the symmetric monoidal structure $\boxtimes$ on {\bf S}--Vect,
see (1.1.3). The identification (1.1.12) shows that any vector space
$X$ gives rise to a tensor functor (``evaluation at $X$")
$$\text{Ev}_X: ({\bold S}\text{--Vect}, \boxtimes)\to (\text{Vect},
\otimes), \quad \Cal V\mapsto \Cal V \langle X\rangle.\leqno (1.6.1)$$
In particular, taking $X=k$, we get a vector space
$$\underline{\Cal V} := \Cal V \langle k\rangle = \bigoplus \Cal V(n)_{
\bold{S}_n}
\leqno (1.6.2)$$
which depends on $\Cal V$ is a multiplicative (with respect to $\boxtimes$)
way. The first part of the following proposition was observed in [F]. 

\proclaim {\quad (1.6.3) Proposition} Let $\Cal P$ be any operad. Then the 
operation
$\boxtimes$ makes $\roman{Mod}$--$\Cal P$, the category of right
$\Cal P$--modules, into a symmetric monoidal abelian category. For any
vector space $X$ the functor $\roman{Ev}_X$ is a tensor functor on
$(\roman{Mod}$--$\Cal P, \boxtimes)$, i.e.,
it takes $\boxtimes$ into $\otimes$.
\endproclaim

 {\bf Proof.} Follows from (1.1.12). 

\medskip

Let now $M$ be a right $\Cal P$--module. We define the operad 
$\text{Op\, End}_{\Cal P}(M)$ as the endomorphism operad of $M$ as an object 
of the symmetric
monoidal category Mod--$\Cal P$:
$$
\text{Op\, End}_{\Cal P}(M)(n) = \text{Hom}_{\text {Mod}-\Cal P}(
M^{\boxtimes n}, M).
\leqno (1.6.4)
$$
Similarly, for two right $\Cal P$--modules $M,N$ we define an
$(\text{Op\, End}_{\Cal P}(N), \text{Op\, End}_{\Cal P}(M))$--bimodule
$\text{Op\, Hom}_{\Cal P}(M,N)$ by
$$\text{Op\,Hom}_{\Cal P}(M,N)(n) = \text{Hom}_{\text {Mod}-\Cal P}
(M^{\boxtimes n}, N). \leqno (1.6.5)$$

\proclaim{\quad (1.6.6) Proposition} We have an isomorphism of operads
$u: \Cal P\to \roman{Op\, End}_{\Cal P}(\Cal P)$. 
\endproclaim

\smallskip

{\bf Proof.} The construction of $u$ is straightforward:
the composition $\mu: \Cal P\circ\Cal P\to\Cal P$
gives, upon restriction to the $n$th summand in (1.1.4), a 
$\bold{S}_n$--equivariant
morphism of right $\Cal P$--modules $\Cal P(n)\otimes_k \Cal P^{\boxtimes n}
\to\Cal P$, and the associativity of $\mu$ implies that $u$ is a morphism
of operads. 
To see that it is injective, look at the action of
$u(p),\,p\in\Cal{P}(n),$ on $1\otimes\dots\otimes 1\in
\Cal{P}(1)^{\otimes n}.$ To check the surjectivity,
consider any element $\varphi =(\varphi_{l_1,\dots ,l_n})
\in \roman { Op\,End}_{\Cal{P}}(\Cal{P})(n).$ Put
$p=\varphi_{1,\dots ,1}(1\otimes\dots \otimes 1)\in\Cal{P}(n)
\subset \Cal P^{\boxtimes n}(n).$
Then $\varphi = u(p).$ In fact, for $q_i\in\Cal{P}(l_i)$
we have
$$
\varphi_{l_1,\dots ,l_n}(q_1\otimes\dots\otimes q_n)=
\varphi_{l_1,\dots ,l_n}(1(q_1)\otimes\dots\otimes 1(q_n))
$$
$$
=(\varphi_{1,\dots ,1}(1\otimes\dots\otimes 1))(q_1,\dots ,q_n)=p(q_1,\dots 
,q_n).
$$
Proposition is proved.

\bigskip

{\bf (1.7) A Lie algebra associated to an operad.}
Let us recall the main ideas of the Tannaka--Krein duality
using [De] as a reference. If $(\Cal C, \otimes)$ is a symmetric
monoidal $k$-linear category and $\Phi: \Cal C\to\text {Vect}$
is a tensor functor, then one can form the group
$\text{Aut}(\Phi)$ of all $\otimes$-automorphisms of $\Phi$.
If $\Cal C$ satisfies additional properties listed in [De]
(i.e., $\Cal C$ is Tannakian, in the terminology of {\it loc. cit.}),
then  $\text{Aut}(\Phi)$  naturally becomes
(the group of $k$-points of) a group
scheme over $k$ and, moreover, $\Cal C$ becomes identified with
the category of regular representations of this group scheme.
Working at the infinitesimal level and weakening
the conditions we have the following 
definition.

\smallskip

\proclaim {\quad (1.7.1) Definition} Let $(\Cal C, \otimes)$ be a
$k$--linear symmetric monoidal category, $\Phi:  \Cal C\to\text {Vect}$
a tensor functor. A derivation of $\Phi$ is a natural transformation
$L: \Phi\to \Phi$ such that for any two objects $A,B\in\Cal C$
the morphism $L_{A\otimes B}: \Phi(A\otimes B)\to \Phi(A\otimes B)$
coincides, after the identification $\Phi(A\otimes B)\simeq \Phi(A)\otimes
\Phi(B)$, with $L_A\otimes 1_{\Phi(B)} + 1_{\Phi(A)} \otimes L_B$. 
\endproclaim

All the derivations of $\Phi$ form, obviously, a Lie algebra which
we denote $\text{Der}(\Phi)$. 

\medskip

We now fix an operad $\Cal P$ and specialize
the above  to the symmetric monoidal category
$(\text{Mod}-\Cal P, \boxtimes)$ and the tensor functor
$\Phi$ given by $\Phi(M)=\underline{M} = M\langle k\rangle$, see Proposition
1.6.3. This category is not Tannakian because it lacks dual
objects, and $\Phi$ is not a fiber functor because
it is not faithful. Nevertheless, Definition 1.7.1 is applicable and gives
some Lie algebra depending only on $\Cal P$.
 We start with a down--to--earth description of this algebra.

\smallskip

For $p\in \Cal P(l)$, $q\in\Cal P(n)$ and $1\leq i\leq l$ set
$$p\circ_i q= p(1, ..., 1, q, 1, ..., 1)\in \Cal P(l+n-1),
\leqno (1.7.2)$$
where $q$ on the right hand side is at the $i$th position, cf. [Mar1].
In a similar way, if $M$ is a right $\Cal P$--module, $m\in M(l)$
and $q\in\Cal P(n)$, we define $m\circ_iq\in M(l+n-1)$, $1\leq i\leq l$.

\proclaim{\quad (1.7.3) Theorem}(a) Let $\Cal P$ be an operad and set for 
$p\in
\Cal P(a+1)$, $q\in\Cal P(b+1)$:
$$p\circ q =  \sum_{i=1}^{a+1} p\circ_i q\quad
\in\quad \Cal P(a+b+1).$$
Then the operation $\circ$ 
induces on the vector space $\underline{\Cal P}$
(see (1.6.2)) the structure of
a pre-Lie algebra in the sense of [Ger] [CL], i.e., it satisfies the following 
property:
$$ (p\circ q)\circ r - p\circ (q\circ r) = (p\circ r)\circ q - p\circ
(r\circ q).$$
In particular, the operation $[p,q] = p\circ q - q\circ p$
makes  $\underline{\Cal P}$ into a Lie algebra which we denote
$\Cal L(\Cal P)$.
 
\smallskip 

(b) For every right $\Cal P$--module $M$ the operation
$$(m,q)\mapsto mq = \sum_{i=1}^{a+1} m\circ_i q \in M(a+b+1),
\quad m\in M(a+1),\, q\in\Cal P(b+1),$$
induces on $\underline M$ the structure of a graded right
module over the pre-Lie algebra $\underline {\Cal P}$ and
hence of a graded right module over the Lie algebra $\Cal L(\Cal 
P)$, i.e., this operation satisfies
the relation of a pre-Lie algebra product except that the first
variable belongs to $\underline {M}$. If $M_1, M_2$
are two right $\Cal P$--modules, then $\underline{M_1\boxtimes M_2}$,
is identified with $\underline{M_1}\otimes \underline{M_2}$ as a
$\Cal L(\Cal P)$--module.

\smallskip 

(c) For any finite--dimensional $\Cal P$--algebra $A$ the Lie algebra
$\Cal L(\Cal P)$ acts on (the affine space underlying) $A$ by polynomial
vector fields. 
\endproclaim

\smallskip

{\bf Proof.} (a) First of all, $\underline{\Cal P}$ is the quotient
of $\bigoplus \Cal P(n)$ by the subspace $J$ spanned by elements of the form
$p(1-s)$, $p\in\Cal{P}(a+1), s\in \bold{S}_{a+1}.$ Using the functional
notation, one easily checks that $J$ is the two--sided ideal
with respect to the product $\circ$. Moreover,  we have, $p\circ_iq\equiv
p\circ_jq\,\text{mod}\,J.$  
\smallskip

The fact that $\circ$ satisfies the identity of a pre-Lie algebra,
follows directly from axioms for the $\circ_i$, given in [Mar1].
Compare with [Ger], where a part of these axioms (not involving
the symmetric group action)  is axiomatized under the name of
a ``pre-Lie system" (\S 5) and the pre-Lie algebra identity is derived
(\S 6).


\smallskip

(b) This statement is checked similarly, and we omit the details.

\smallskip

(c) Notice that $A^\star = \{(A^*)^{\otimes n}\}_{n\geq 0}$ 
is naturally a right $\Cal P$--module, and $\underline{A^\star}=
S(A^*)$ is the algebra of polynomial functions on (the affine space
underlying) $A$. So $S(A^*)$ becomes a $\Cal L(\Cal P)$-module, and to
prove our assertion we need only to check that  
$\Cal L(\Cal P)$ acts by algebra derivations. But this follows
from the compatibility with tensor products. 

\medskip

{\bf (1.7.4) Remark.} The proof above
shows that the direct sum $\Lambda(\Cal P) =\bigoplus \Cal P(a+1)$
also forms a graded (pre-)Lie algebra, of which $\Cal L({\Cal P})$
is a quotient. The Gerstenhaber bracket on the Hochschild
complex of an associative algebra $A$, see [Ger], is a particular case of the
construction of  $\Lambda(\Cal P)$  but for 
the operad
 $\text{Op \, End}(A[-1])$ in the category of graded vector spaces.

\smallskip 

\proclaim{\quad (1.7.5) Theorem} The  Lie algebra $\Cal L(\Cal P)$
is identified with $\roman{Der}(\Phi)$ for the tensor functor described
above.
\endproclaim

{\bf Proof.} Part (b) of Theorem 1.6.3 means that
we have a morphism of Lie algebras $u: \Cal L(\Cal P)\to
\text{Der}(\Phi)$. To see that $u$ is injective it is
enough to consider an action of $p\in \Cal P(a+1)$ on $\underline{\Cal P}$
where $\Cal P$ is
considered as a module over itself. It remains to prove surjectivity,
i.e., that any derivation $D: \Phi\to \Phi$ comes from some element 
of $\Cal L(\Cal P)$. To see this, consider the $D$--action on
$\underline{\Cal P}$ and let $q = D(1)$ be the image of $1\in \Cal P(1)$.
We claim that $D=u(q)$. Indeed, let $M$ be an arbitrary right
$\Cal P$--module and $m\in M(a+1)$. Then the associativity of the right $\Cal 
P$--action
on $M$ can be expressed as follows.
Consider the natural morphism
$$M(a+1)\otimes_k \Cal P^{\boxtimes (a+1)}\to M\leqno (1.7.6)$$
given by restricting the action $M\circ \Cal P\to M$, see (1.1.4).
Let us view the source of this morphism as a right $\Cal P$--module by 
considering
$M(n)$ as just the vector space of multiplicities of the $\Cal P$--module
$\Cal P^{\boxtimes (a+1)}$.
Then (1.7.6) is a morphism of modules.

\smallskip

 Now, because $D$ is a derivation,
applying $D$ to $m\otimes 1\otimes... \otimes 1$
in the source of (1.7.6) we get
$$\sum_{i=1}^{a+1} m\otimes 1\otimes ... \otimes 1\otimes
 q\otimes 1\otimes ... \otimes 1$$
with $q$ at the $i$--th place in the $i$--th term. Since (1.7.6) is a morphism 
of
modules, we find that the $D$--action on $m$ is given by the formula
in (1.7.3)(b). Theorem is proved.

\medskip

{\bf (1.7.7) Examples.} (a) Let $V$ be a finite-dimensional vector
space and $\Cal P=\text{Op\, End}(V)$. Then
$$
\underline{\Cal P} = \bigoplus_{n\geq 0} S^n(V^*)\otimes V
$$
can be naturally identified with $\text{Der}\,k[V^*]$, the space of
polynomial vector fields on $V$. Moreover, the bracket $[p,q]$
is identified with the usual Lie bracket of vector fields, so
$\Cal L(\text{Op\, End}(V))\simeq \text{Der}\,k[V^*]$ as a Lie algebra,
the action (1.7.3)(c) being the tautological one.

The fact that the Lie algebra of vector fields on $V$ comes in
fact from a pre-Lie algebra, was pointed out in [Ku].

\smallskip

(b) Taking $V=k$ in (a), we get $\Cal P = \Cal Com$. Thus,
$\Cal L(\Cal Com) = \text{Der}\,k[x]$ is
 Lie algebra of polynomial vector fields on the line. The commutation
relations for the canonical generators
$e_n\in \Cal Com(n)$ is directly found to be
$$[e_m, e_n] = (m-n) e_{m+n-1}.$$
 In other words, $e_n$  
corresponds to the vector field $x^n(d/dx)$.

Further, take $\Cal P=\Cal Ass$ to be the associative operad. Then
$\Cal P(n)_{\bold{S}_n} = k$ and we find that $\Cal L(\Cal Ass)$
is isomorphic to $\Cal L(\Cal Com)$, i.e., to $\text{Der}\,k[x]$. 
As a corollary, we get part (a) of the following fact.

\smallskip

\proclaim {\quad (1.7.8) Theorem} (a) Let $A$ be a
 finite--dimensional associative
algebra. Then the Lie algebra 
 $\text{Der}\,k[x]$ acts by
polynomial vector fields on $A$. Explicitly, $e_n=x^n (d/dx)$ acts
by the vector field
$$E_n = \biggl( X^n, \, {\partial\over\partial X}\biggr).$$
Here $X^n$ is the $A$--valued polynomial function on $A$ which raises
every element to the $n$--th power, ${\partial\over\partial X}$ is the 
canonical
$A^*$--valued vector field on $A$ (corresponding to
 $\text{Id}\in A^*\otimes_k A$) and $(-,-)$ is the structure pairing between
$A$ and $A^*$. 

\smallskip

(b) Let $A^\times \subset A$ be the set of invertible elements.
Let $\operatorname{Der}\, k[x, x^{-1}]$ be the Witt algebra of regular
vector fields on $k^\times = k-\{0\}$, with the basis $e_n=x^n(d/dx),
n\in {\bold Z}$. Then the operators $E_n$ defined similarly to (a),
give an action of   $\operatorname{Der}\, k[x, x^{-1}]$ by regular
vector fields on $A^\times$. 

\endproclaim

This fact for $A=\text{Mat}\,(r, \bold{C})$ was noticed by N.~Wallach ([W]).
The proof of part (b) in general is left to the reader.

\medskip

As we saw in Proposition 1.3.8, any augmented
simplicial vector space gives rise to a right
$\Cal Ass$--module and thus, by the above, it defines a representation of
$\text{Der}\,k[x]$. Let us describe this action explicitly. Let
$V_\bullet$ be an augmented simplicial vector space, with face and degeneracy operators
$\partial_{n,i},
s_{n,i},$ see (1.3.7). 
Let $C(V_\bullet) = \bigoplus_n V_n$ be the graded vector space
associated to $V$. Define the linear operators $L_p: C(V_\bullet)\to
C(V_\bullet)$ of degree $p-1$, $p\geq 0$, as follows. On the summand $V_n
\subset C(V_\bullet)$ we set
$$L_0=\sum_{i=0}^n \partial_{n,i}, \quad L_1 = (n+1)\cdot \text{Id},
\quad L_p = \sum_{i=0}^n s_{n-p+1, i} \dots s_{n+1, i} s_{n,i}, \, \,\, p\geq 
1.$$

\proclaim {\quad (1.7.9) Theorem} The operators $L_p$ satisfy the commutation
relations $[L_p, L_q] = (p-q)L_{p+q-1}$ and thus make $C(V_\bullet)$
into a $\text{Der}\,k[x]$--module. 
\endproclaim

The proof can be deduced directly from the standard simplicial
identities between faces and degenerations. 

\bigskip

{\bf (1.7.10) Example.} Here we describe the Lie algebra
$\Cal{L}(H_*\overline{M})$ for the operad $H_*\overline{M}$
introduced in (1.3.9).
Consider a stable tree with $n+1$ tails (flags
which are not halves of the edges).
Stability means that
there are at least three flags at each vertex.
Any numbering of flags by $[n]$ 
determines an irreducible closed submanifold of $\overline{M}_{0,[n]}$
whose open dense subset parametrizes curves of genus zero
of the respective combinatorial type. The image of the homology
class of this submanifold in $\Cal{L}(H_*\overline{M})$
depends only on the isomorphism class of the respective rooted
tree where root is the flag labeled by $0.$
For any such (isomorphism class) $\tau$ denote by $L(\tau )$
the corresponding element of the Lie algebra. Then we have:

\smallskip

(i) $\Cal{L}(H_*\overline{M})$ is spanned by all $L(\tau ).$
Besides the general $\bold{Z}$--grading (see (1.6.2)), it has
an additional $\bold{Z}$--grading by algebraic dimension of $L(\tau )$
which is $n-2$ less the number of edges of $\tau$.

\smallskip

(ii) The bracket is defined by
$$
[L(\sigma ), L(\tau )] = \sum_i L(\sigma\circ_i\tau )-
\sum_j L(\tau\circ_j\sigma ).
\leqno (1.7.11)
$$
Here $i$ (resp. $j$) runs over non--root flags of $\sigma$
(resp. of $\tau$), and $\sigma\circ_i\tau$ denotes the tree obtained
from $\sigma$ by gluing its flag $i$ to the root of $\tau$.

\smallskip

(iii) The elements $L(\tau )$ are lineraly independent.
This follows from the fact that all linear relations
between $L(\tau )$ can be obtained by symmetrizing
the linear relations between the boundary 
homology classes described in [Man2], sec. III.4.7,
and this symmetrization produces zero.

\smallskip

All of this in the final count follows from Keel's 
description of the Chow and homology groups
of $\overline{M}_{0,[n]}$ ([Ke]).  

\medskip

For a general $H_*\overline{M}$--algebra $H$ there is nothing to
add to the description of vector fields given in the proof of
(1.7.3) (c). Let us therefore consider the case of the cyclic
algebra $H$, for example, quantum cohomology.
In this case $H$ comes with a richer structure, namely it is endowed with a 
non--degenerate
symmetric pairing $g$ (Poincar\'e form), and the algebra structure 
geometrically
appears in the guise of Gromov--Witten $\bold{S}_{n+1}$--equivariant
maps 
$$
I_{n+1}:\, H^{\otimes (n+1)}\to H^*(\overline{M}_{0,[n]},k).
$$
The operadic action needed to define $L(\tau )$
$$
\mu_n :\, H_*(\overline{M}_{0,[n]})\otimes H^{\otimes n}\to H
$$
is obtained from $I_{n+1}$ by partial dualization with the help of $g.$

\bigskip

{\bf (1.8) The role of the product $\otimes$.} Consider now the
monoidal structure $\otimes$ on {\bf S}--Vect, see (1.1.2).

\smallskip 
 
If $\Cal P_1, \Cal P_2$ are operads, then the distributivity maps
(1.1.13) make $\Cal P_1\otimes \Cal P_2$ into an operad. If $A_i$ is a
$\Cal P_i$--algebra, $i=1,2$, then $A_1\otimes_k A_2$ is a $\Cal P_1
\otimes\Cal P_2$--algebra via (1.1.14). Similarly, 
if $M_i$ is a left (resp. right) $\Cal P_i$--module, i=1,2, then
$M_1\otimes M_2$ is  a left (resp. right) $\Cal P_1\otimes \Cal P_2$--module.

\smallskip

Let $\Cal P$ be an operad, $M$ be a right $\Cal P$--module
and $V$ be a vector space.  Then we have the right $\Cal P$--module
$V\otimes M$, see (1.3.3). Notice that if $N, W$ are
another right $\Cal P$--module and a vector space, then we have
an identification
$$(V\otimes M)\boxtimes(W\otimes N) = (V\otimes_k W)\otimes(M\boxtimes N),
\leqno (1.8.1)$$
which implies the following. 

\smallskip

\proclaim{\quad (1.8.2) Proposition} We have an isomorphism of operads
$$\roman{Op\, End}_{\Cal P}(V\otimes M) \simeq \roman{Op\, End}(V)\otimes
\roman{Op\, End}_{\Cal P}(M)$$ 
and an isomorphism of bimodules
$$\roman{Op\, Hom}_{\Cal P}(V\otimes M, W\otimes N) \simeq
\roman{Op\, Hom}(V,W)\otimes\roman{Op\, Hom}_{\Cal P}(M,N).$$
\endproclaim

\proclaim{\quad (1.8.3) Definition} The $d$ by $d$ matrix operad over $\Cal P$
is defined by
$$\roman { Mat}\,(d,\Cal{P}) = \Cal{P}\otimes \roman { Op\,End}\,(k^d) = 
\roman{Op\, End}_{\Cal P}(\Cal P^{\oplus d})$$
(the last equality being  a consequence of (1.8.2) and Proposition 1.6.6).
\endproclaim

Next, let $\Cal P_i, i=1,2$ be operads, $M_i$ be a right $\Cal P_i$--module
and $N_i$ be a left $\Cal P_i$--module. Then the distributivity maps
(1.1.12) give rise to a morphism of {\bf S}--modules
$$(M_1\otimes M_2)\circ_{\Cal P_1\otimes\Cal P_2}(N_1\otimes N_2)
\to (M_1\circ_{\Cal P_1}N_1)\otimes (M_2\circ_{\Cal P_2} N_2).\leqno(1.8.4)$$
In general this need not be an isomoprhism. We will be interested in
the case when $\Cal P$ is a fixed operad, $M_1=N_1=\Cal P_1=\Cal P$,
while $\Cal P_2=\text{Op\, End}(V)$, $\dim\,V <\infty$, $M_2=V^\star$ and 
$N_2=V\otimes\Cal Com$. Notice that in this case $M_2\circ_{\Cal P_2} N_2
\simeq\Cal Com$, see Example 1.5.5(a). Let us denote $V\otimes\Cal Com$
by $\widetilde V$. The specialization of the morphism
(1.8.4) to our case is then:
$$u: (\Cal P\otimes V^\star)\circ_{\Cal P\otimes \roman{Op\, End}(V)}
(\Cal P\otimes\widetilde V)\to \Cal P\otimes(V^\star\circ_{\roman{Op\, End}(V)}
\widetilde V) = \Cal P.
\leqno (1.8.5)$$

\smallskip

\proclaim {\quad(1.8.6) Proposition} The morphism $u$ in (1.8.5) is an
isomorphism of {\bf S}--modules.

\endproclaim

\smallskip

{\bf Proof.} Let $\dim(V)=d$. It is enough to show that
$u^{\oplus d}$, the direct sum of $d$ copies of $u$, is an isomorphism.
On the other hand, $(V^\star)^{\oplus d}$ is identified, as a right
$\text{Op\, End}(V)$--module, with $\text{Op\, End}(V)$ itself and hence
$(\Cal P\otimes V^\star)^{\oplus d}$ is identified, as a 
$\Cal P\otimes \text{Op\, End}(V)$--module, with 
$\Cal P\otimes \text{Op\, End}(V)$ itself. It follows that $u^{\oplus d}$
is identified with the morphism
$$\Cal P\otimes \text{Op\, End}(V)\circ_{\Cal P\otimes \text{Op\, End}(V)}
(\Cal P\otimes\widetilde V) \to \Cal P\otimes \left(\text{Op\, End}(V)
\circ_{\text{Op\, End}(V)} \widetilde V\right)$$
which is just the identity map $\Cal P\otimes\widetilde V\to
\Cal P\otimes\widetilde V$. Thus $u^{\oplus d}$ is an isomorphism and hence
$u$ is an isomorphism. 

\bigskip

{\bf (1.9) Morita equivalence.}
 Consider an operad $\Cal{P}$ and a
$\Cal{P}$--module $M$. Put $\Cal{Q}=\roman {Op\,End}_{\Cal{P}}(M).$
Then $M$ is a $(\Cal{Q},\Cal{P})$--bimodule,
so that the construction of (1.4) gives us a functor
$f_M:\,\Cal{P}$--Alg$\to \Cal{Q}$--Alg.

\medskip

\proclaim{\quad (1.9.1) Theorem} If $M$ is free as $\Cal{P}$--module,
that is, isomorphic to $\Cal{P}^{\oplus d},$ then
$f_M$ is an equivalence of categories.
\endproclaim

\smallskip

{\bf Proof.}  As above, put $V=k^d$ so that $M=\Cal{P}\otimes_k V  = 
 \Cal P\otimes\widetilde V.$ Let $M^* = \text{Op\, Hom}_{\Cal P}(M, \Cal P) = 
\Cal P\otimes V^*$, the last identification following from Proposition 1.8.2. 
Then $M^*$ is a $(\Cal P, \Cal Q)$--bimodule and 
  we have the functor
$f_{M^*}:\,\Cal{Q}$--Alg$\to \Cal{P}$--Alg.
We claim that $f_M$ and $f_{M^*}$ are mutually inverse.

\smallskip

To prove this, it suffices to construct the isomorphisms
$$
M^*\circ_{\Cal{Q}}M\cong \Cal{P}\quad (\roman{as}\
(\Cal{P},\Cal{P})\roman{-bimodules}).
\leqno(1.9.2)
$$
$$
M\circ_{\Cal{P}}M^*\cong \Cal{Q}\quad (\roman{as}\
(\Cal{Q},\Cal{Q})\roman{-bimodules}).
\leqno(1.9.3)
$$

Now, (1.9.2) is the content of Proposition 1.8.6. As for (1.9.3),
we have, first of all, a morphism
$$\psi: M\circ_{\Cal P} M^* = (\Cal P\otimes\widetilde V)\circ_{\Cal P\otimes
\Cal Com}(\Cal P\otimes V^\star) \to \Cal P\otimes(V\circ_{\Cal Com} V^\star),
\leqno (1.9.4)$$
a particular case of (1.8.4). Note that the natural morphism of
$\text{Op\, End}(V)$--bimodules 
 $$\phi: \widetilde V\circ_{\Cal Com}V^\star\to \text{Op\, End}\,(V)$$
is an isomorphism. This is because $\widetilde V$ is a free right
$\Cal Com$--module, so $\phi$ is isomorphic to the direct sum
of $d=\dim(V)$ copies of the isomorphism $\Cal Com\circ_{\Cal Com}V^\star\to
V^\star$. Thus (1.9.4) can be regarded as a morphism
$$\psi: M\circ_{\Cal P}M^*\to \Cal P\otimes\text{Op\, End}(V) = \Cal Q.$$
To see that this is an isomorphism, it is enough to notice that
$\psi$ is isomorphic  (as a morphism of {\bf S}--modules) to the direct sum of
$d$ copies of the morphism
$$\Cal P\circ_{\Cal P}(\Cal P\otimes V^\star) \to \Cal P\otimes V^\star,$$
which is an isomorphism. Theorem 1.9.1 is proved. 

\bigskip

{\bf (1.10) The super--version.} Most of the above
constructions can be performed in any $k$--linear
abelian symmetric monoidal category $\Cal C$, instead of Vect.
We will be particularly interested in the category
SVect of super-vector spaces. Recall that objects of
SVect are $\bold Z/2$--graded $k$--vector spaces,
the tensor product is the usual graded one and the symmetry isomorphism
$V\otimes W\to W\otimes V$ involves the Koszul sign.

\smallskip 

If $\Cal P$ is an operad in SVect, then $\Cal P$--algebras
(in SVect) will be called $\Cal P$--superalgebras and their
category will be denoted $\Cal P$--SAlg.

\smallskip

All the statements given above for $k$--linear
operads generalize to SVect without difficulty.
In particular, for any super--vector space $V$ we have the endomorphism
operad $\text{Op\, End}(V)$ in SVect, for any two super--vector
spaces $V,W$ we have a bimodule $\text{Op\, Hom}(V, W)$ and so on.

\smallskip

We will need a superversion of the Morita theorem and of a more precise
statement underlying it. Let $\Cal P$ be any operad in SVect,
$M$ any right $\Cal P$--module, $\Cal Q = \text{Op\, End}_{\Cal P}(M)$
and $M^* = \text{Op\, Hom}_{\Cal P}(M, \Cal P)$. Then, as before,
we have a  morphism
$$M^*\circ_{\Cal Q} M\to \Cal P \leqno (1.10.1)$$
of $(\Cal P,\Cal P)$-bimodules and a morphism
$$M\circ_{\Cal P} M^*\to \Cal Q \leqno (1.10.2)$$
of $(\Cal Q, \Cal Q)$-bimodules.

\smallskip 

\proclaim{\quad (1.10.3) Theorem} If $M$ is free, i.e., isomorphic to
$\Cal P\otimes_k V$ where $V$ is a finite--dimensional supervector
space, then the morphisms (1.10.1--2) are isomorphisms and hence
the functor $f_M$ is an equivalence between $\Cal P$--SAlg and
$\Cal Q$-SAlg. 

\endproclaim

The proof is obtained by performing the same steps as for Theorem
1.9.1, with easy modifications to accomodate the $\bold Z/2$-grading.

\newpage

\centerline {\bf \S 2. Differential algebras.}

\bigskip

 {\bf (2.1) Operadic approach to nonlinear differential equations.}
Let $X$ be a smooth complex algebraic variety of dimension $p$. Denote by
${\Cal D}(n) = {\Cal D}_X(n)$ the sheaf of $n$--linear multidifferential
operators 
$$(u_1, ..., u_n) \mapsto L(u_1, ..., u_n), \quad u_i\in {\Cal O}_X$$
acting on regular functions on $X$. If $(x_1, ..., x_p)$ 
is a local coordinate system on
$X$, and $\partial^I = \partial_{x_1}^{i_1}...\partial_{x_p}^{i_p}$,
$I=(i_1, ..., i_p),$ is the iterated partial derivative
corresponding to a multiindex $I$, then a local section $L$ of
${\Cal D}(n)$ acts on functions as follows
$$L(u_1, ..., u_n) = \sum_{I(1), ..., I(n)} f_{I(1), ..., I(n)}(x)
(\partial^{I(1)} u_1)\cdot ... \cdot(\partial^{I(n)} u_n), 
\quad  f_{I(1), ..., I(n)}(x)\in \Cal O_X.$$
Thus $\Cal D(1)$ is the usual sheaf of rings of linear differential operators
on $X$. The collection $(\Cal D_X(n))_{n\geq 0}$ forms a sheaf of
operads on $X$, with the composition given by
superposition of multilinear differential operators. 
As with any operad, $\Cal D_X(n)$ is endowed with a left $\Cal D_X(1)$--module 
(in particular, a $\Cal O_X$-module) structure
and with $n$ commuting right $\Cal D_X(1)$--module (in particular, 
$\Cal O_X$--module) structures. 

\smallskip
 
 As is well known, the study of sheaves of modules over
 the sheaf of rings $\Cal D_X(1)$ provides a fruitful algebraic
approach to the
 theory of linear differential equations. The introduction of
the sheaf of operads $\Cal D_X$  and its sheaves of algebras
provides a similar algebraic language for systems of
nonlinear differential equations. More precisely, let us give
the following definition.

\smallskip

\proclaim {\quad \bf (2.1.1) Definition} A system of  
differential equations on $X$ 
(polynomial in the derivatives) is a sheaf $\Cal A$ of $\Cal D_X$--algebras
which is locally given by finitely many generators and relations. 
A  solution to a system of equations given by $\Cal A$ is a morphism
of $\Cal D_X$--algebras ${\Cal A}\to\Cal O_{X}$ (note that $\Cal O_{X}$
 is
a $\Cal D_X$--algebra by definition).
\endproclaim

To understand how this is related with the intuitive concept of a system
of differential equations, consider the case of
one unknown function $u$ and a system of
differential equations
$$P_\nu\biggl( u, {\partial u\over\partial x_i}, {\partial^2 u\over
\partial x_i\partial x_j}, \,\,\, ...\biggr)=0\leqno (2.1.2)$$
where $P_\nu$ is a polynomial in $u$ and its derivatives
with coefficients in $\Cal O_X$. We polarize each homogeneous
component of $P_\nu$ to a multilinear differential operator
thus writing each equation from (2.1.2) as
$$\sum_n L_{\nu, n}(u, ..., u)=0.\leqno (2.1.3)$$
where each $L_{\nu, n}$ is an $n$--ary multilinear differential operator,
i.e., a section of $\Cal D_X(n)$.  Then we form a $\Cal D_X$--algebra $\Cal A$
on one generator $U$ subject only to the relations (2.1.3). A $\Cal 
D_X$--algebra
homomorphism $\Cal A\to \Cal O_X$ is clearly the same as a solution to
(2.1.2), as $U$ should go to some function $u\in \Cal O_X$ so that
the defining relations of $U$ in $\Cal A$ are satisfied in $\Cal O_X$. 
Similarly, a system of equations on $n$ unknown functions
is translated into an algebra with $n$ generators.

\smallskip

The next proposition can be checked by straightforward
manipulations in local coordinates.

\smallskip

\proclaim {\quad (2.1.4) Proposition} (a)  The filtration $F_\bullet \Cal 
D_X(n)$
by total number of derivatives is compatible with the operad structure
(i.e., makes $\Cal D_X$ into an operad in the category
of filtered quasicoherent sheaves of $\Cal O_X$-modules).
In particular, $F_0\Cal D_X$ is a sheaf of operads which can be identified
with $\Cal Com\otimes_{\bold C}{\Cal O}_X$, and $\roman {gr}_\bullet^F
\Cal D_X$ is an operad in the category of graded sheaves.

\smallskip

(b) A $\Cal D_X$-algebra is the same as a commutative $\Cal O_X$-algebra
(via $\Cal Com\otimes_{\bold C}{\Cal O}_X = F_0\Cal D_X\i \Cal D_X$)
which is made into a left $\Cal D_X(1)$-module such that vector
fields in $\Cal D_X(1)$ act by algebra derivations.
\endproclaim

Part (b) shows that our approach is in fact identical to the 
classical approach of ``differential algebra" of Ritt [Ri].
However, the operadic point of view seems to present
several advantages. Let us explain, for example, the analog
of the fundamental fact that
  the algebra
$\roman {gr}_\bullet^F \Cal D_X(1)$ is commutative
and can be identified with the algebra of functions on $T^*X$.
Recall that an algebra $A$ is commutative if and only
if the multiplication $A\otimes A\to A$ is an algebra morphism,
i.e., $A$ is an algebra in the category of algebras.

\proclaim {\quad \bf (2.1.5) Proposition} (a) Each sheaf
 $\roman{gr}_\bullet^F \Cal D_X(n)$
has a natural structure of a commutative algebra so that
$\roman{gr}_\bullet^F \Cal D_X$ becomes an operad in the category of
algebras.

\smallskip

(b) The spectrum of the algebra $\roman{gr}_\bullet^F \Cal D_X(n)$
is identified with $(T^*_X)^{\oplus n}$.

\smallskip

(c) The operad structure on $\roman{gr}_\bullet^F \Cal D_X$
is induced by a structure of the cooperad (in the category of
algebraic varieties) on the collection of $(T^*_X)^{\oplus n}$ given
by the maps
$$\nu_{a_1, ..., a_n}: (T^*_X)^{\oplus (a_1+...+a_n)}
\to (T^*_X)^{\oplus n}\times (T^*_X)^{\oplus a_1}\times ...
\times (T^*_X)^{\oplus a_n},\leqno (2.1.6)$$
$$(\xi_1, ..., \xi_{a_1+...+a_n})\mapsto $$
$$ \left(
\biggl(\sum_{i=1}^{a_1}\xi_i, \sum_{i=a_1+1}^{a_1+a_2}\xi_i, ..., 
\sum_{i=a_1+...+a_{n-1}+1}^{a_1+...+a_n}\xi_i\biggr), \,\,
(\xi_1, ..., \xi_{a_1}), \,\, ..., \,\,
(\xi_{a_1+...+a_{n-1}+1}, ..., \xi_{a_1+...+a_n})\right).$$
\endproclaim

Again, the proof is straightforward. 

\medskip

 {\bf (2.1.7) Remarks.}
(a) Note that the  construction (2.1.6) can be defined in fact for any abelian
 group
$A$ (instead of $T^*_X$), making $(A^{\oplus n})_{n\geq 0}$ into a cooperad.

\smallskip

(b) Proposition 2.1.5(c) suggests a natural algebraic approach to
microlocalization of nonlinear equations: instead of an open
cone in $T^*_X$ in linear theory, the microlocalization should
be done with respect to  a sub-cooperad
in  $\{(T^*_X)^{\oplus n}\}$.

\medskip

 {\bf (2.2) $\Cal D_X$ as an endomorphism operad.}
We can view the operad $\Cal D_X$, though not formally as  particular
case, but as a close analog of the construction of the endomorphism
operad $\roman{Op\, End}\,(V)$ from (1.2.6): we take for $V$ the
space (sheaf) $\Cal O_X$ of functions on $X$ and consider
only local endomorphisms. All the abstract constructions of \S 1
can be performed in this situation and have interesting meaning.

\smallskip

First of all, the role of the dual space $V^*$ is played by
$\omega_X$, the sheaf of volume forms on $X$, and the pairing
between $V$ and $V^*$ is replaced by the Serre duality. The analog of
the fact that $\roman{Op\,End}\,(V)(n) = (V^*)^{\otimes n}\otimes V$
is given by the following coordinate-free definition of  
  $\Cal D_X(n)$
in terms of local cohomology sheaves, generalizing Sato's definition of 
linear differential operators [SKK].

\smallskip

\proclaim {\quad (2.2.1) Proposition}
Consider the
$(n+1)$-fold Cartesian product $X^{n+1}$ whose factors will be labeled
by integers $0,1,..., n$ and let $p_i: X^{n+1}\to X$,
$i=0,...,n$, be the projections. Let also $\Delta\i X^{n+1}$
be the image of the diagonal embedding of $X$. Then we have a natural
identification
$$\Cal D_X(n) = \underline{H}^{n\dim(X)}_\Delta\biggl(
\bigotimes_{i=1}^n p_i^*\omega_X \biggr)$$
\endproclaim

\smallskip

 {\bf Proof.} We separate the question in two. First, let $f: X\to Y$
be any morphism of smooth algebraic varieties. One has then the sheaf
$$\Cal D_{X\to Y} = \Cal O_X\otimes_{f^{-1}\Cal O_Y} f^{-1}\Cal D_Y$$
of $(\Cal D_X(1), f^{-1}\Cal D_Y(1))$-bimodules on $X$, 
see, e.g., [Bo] IV, \S 4.2
or [Ka] p. 24. Its sections can be seen as linear
differential operators
 $f^{-1}\Cal O_Y\to \Cal O_X$, see [Ka], {\it loc. cit.} 
Our proposition is then implied by the next two facts.

\proclaim{(2.2.2) Lemma}
If $\Gamma(f)\subset X\times Y$ is the graph of $f$ and $p_Y:
X\times Y\to Y$ is the projection, then $D_{X\to Y}$
is naturally identified with 
$\underline{H}^{\dim(Y)}_{\Gamma(f)}(p_Y^*\omega_Y)$.

\endproclaim

\proclaim{(2.2.3) Lemma} If $Y=X^n$ 
and $f: X\to X^n$ is the diagonal embedding,
then $D_{X\to X^n} = \Cal D_X(n)$.

\endproclaim

Lemma 2.2.2 is fairly classical, see, e.g., [Ka], p.24 or (for a similar
fact about pseudodifferential operators) [SKK], p. 329. 

To see Lemma 2.2.3, notice that for $n$ functions $u_1, ..., u_n
\in\Gamma(U, \Cal O_X)$, we have the function $u_1\otimes ... \otimes
u_n\in\Gamma(U^n, \Cal O_{X^n})$, which then gives a section
of $f^{-1}\Cal O_{X^n}$ over $U$, still denoted  $u_1\otimes ... \otimes
u_n$. Now, applying linear differential operators $f^{-1}\Cal O_{X^n}\to
\Cal O_X$ (i.e., sections of $\Cal D_{X\to X^n}$) to sections of
the form  $u_1\otimes ... \otimes
u_n$, we get precisely all the $n$-linear differential operators
$(\Cal O_X)^n\to\Cal O_X$, as one can easily see in local
coordinates. Lemma 2.2.3 and Proposition 2.2.1 are proved.

\medskip

As we are almost in the endomorphism operad situation,
let us further compare it to the framework of (1.3.1).
Writing $\sim$ to denote the analogous objects, we have
$\Cal P\sim\Cal D_X$, $V\sim \Cal O_X$. Next, the analog of
the module $V^\star$ from (1.3.1) looks as follows. Consider $X^n$, with
projections $p_i: X^n\to X$, $i=1, ..., n$ and the
diagonal $\Delta \subset X^n$, so that $\Delta$ is identified with $X$.
Let
$$\Omega(n) = \underline{H}_\Delta^{(n-1)\dim(X)}\biggl(
\bigotimes_{i=1}^n p_i^*\omega_X\biggr) = \omega_X\otimes_{\Cal O_X}
\Cal D_X(n-1). \leqno (2.2.4)$$

\smallskip

\proclaim {\quad (2.2.5) Proposition} The collection
 $\Omega = \{\Omega(n)\}_{n\geq 1}$ is naturally a right $\Cal D_X$-module
in the sense of (1.2). 
\endproclaim

\smallskip

The proof is straightforward.

\smallskip

Proposition 2.2.5
generalizes the fact that $\omega_X$ is a right module over the
algebra $\Cal D_X(1)$. 

\medskip

\medskip

 {\bf (2.3)  $\Cal D$-algebras on supermanifolds.}
Let now $X$ be a supermanifold of dimension $(p|q)$, see [Man1]. Ch.~4.
Thus $X$ consists of an ordinary smooth complex algebraic variety
$X_{\roman {red}}$ of dimension $p$ and a sheaf $\Cal O_X$ of
${\bold Z}/2$-graded supercommutative algebras which
is locally isomorphic to $\Cal O_{X_{\roman {red}}}\otimes_{\bold C}
\Lambda[\xi_1, \dots , \xi_q]$. Here $\Lambda[\xi_1, \dots , \xi_q]$
is the exterior algebra on the generators $\xi_i$.
By $\omega_X$ we will understand the sheaf of volume forms
in the sense of Berezin, see {\it loc.~cit}.
All the constructions  and statements of (2.1)--(2.2) can be immediately 
generalized to
the super case, giving a sheaf of operads (in SVect) $\Cal D_X$ on $X$
  of which
$\Cal O_X$ is a sheaf of superalgebras. 
 
In particular, consider the case $p=0$, i.e., $X={\bold C}^{(0|q)} = 
\roman {Spec}\, \Lambda$, where $\Lambda=\Lambda[\xi_1, ..., \xi_q]$.
Note that $\Lambda$ is finite--dimensional super--vector space
over {\bf C}, of dimension
$(2^{q-1}| 2^{q-1})$ and $X_{\roman {red}}$ is in this case just a point, so a 
sheaf on
it is just a vector space.

\proclaim{\quad  (2.3.1) Proposition} (a) For $X={\bold C}^{(0|q)}$ the sheaf
of operads $\Cal D_X$ on $X_{\roman {red}} = \{\roman  pt\}$ is
 $\roman {Op\,End}(\Lambda)$, the endomorphism operad
of $\Lambda$ considered as a super-vector space.

\smallskip

(b) If $X=X_{\roman {red}}\times {\bold C}^{(0|q)}$ is a split
supermanifold, then $\Cal D_X \simeq \roman  {\Cal D}_{X_{\roman 
{red}}}\otimes_{\bold C}\roman {Op\,End}(\Lambda)$ is
isomorphic to a matrix operad over $\Cal D_{X_{\roman {red}}}$.

 \endproclaim

 {\bf Proof.} Part (b) follows from (a). To see (a), 
notice first that $\Cal D$-action on $\Cal O$, i.e., on $\Lambda$,
gives an operad morphism $\Cal D\mapsto \text{Op\, End}(\Lambda)$. 
Let $\xi^I$, $I=(1\leq i_1 < ... < i_m \leq q)$, be the monomial
basis in $\Lambda$. Let $F_{I(1), ..., I(n), J}\in\text{Hom}_{\bold 
C}(\Lambda^{\otimes n},
\Lambda)$ be the ``matrix unit" which takes $\xi^{I(1)}\otimes ... \otimes
\xi^{I(n)}$ to $\xi^J$ and all the other tensor products
of monomials to 0. Such operators  form a basis
in $\text{Hom}_{\bold C}(\Lambda^{\otimes n},
\Lambda)$. On the other hand, the space $\Cal D(n)$ of all 
$n$--linear differential operators on $\Lambda$ has a basis formed by operators
$$L_{I(1), ..., I(n), J}, \quad u\mapsto 
 \xi^J (\partial^{I(1)}u) ... (\partial^{I(n)}u).$$
Consider the filtration of $\Lambda$ by the powers of the
maximal ideal $(\xi_1, ..., \xi_q)$ and the induced filtration
on $\text{Hom}_{\bold C}(\Lambda^{\otimes n}, \Lambda)$.
Modulo this filtration, the action of $L_{I(1), ..., I(n), J}$
is given by $F_{I(1), ..., I(n), J}$, whence the statement.

\medskip

Theorem 1.10.3 is now formally applicable in the situation of (2.3.1)(b)
but we prefer to make a slightly more general statement.
Let $\Cal D_X$--SAlg be the category of sheaves of $\Cal D_X$--superalgebras,
see (1.10).
  Note that  even if $X$ happens to be purely
even, then a $\Cal D_X$-superalgebra is still required to be
${\bold Z}/2$-graded. 

\smallskip

\proclaim{\quad (2.3.2) Theorem} The categories $\Cal D_X$--SAlg and
$\Cal D_{X_{\roman{red}}}$--SAlg are equivalent.
\endproclaim

\smallskip

This is an operadic extension of the result of I. Penkov [P] on
$\Cal D_X(1)$-modules. 

\medskip

 {\bf Proof.}
Consider the natural embedding of the supermanifolds $i: X_{\roman{red}}
\to X$. On the level of underlying spaces this is the identity
and the source and target of $i$ differ only by the sheaves of
rings: $\Cal O_{X_{\roman{red}}}$ as opposed to $\Cal O_X$.
Accordingly, both $\Cal D_X$ and ${\Cal D}_{X_{\roman {red}}}$ are sheaves
of operads in SVect on the same underlying space $X_{\roman{red}}$. 
This means that we can apply the formalism of $\text{Op\, Hom}$
bimodules of Section 1, if we understand them in a sheaf--theoretic sense.
 Let 
$$\Cal D_{\to} = \{ D_{\to}(n)\} = \underline{\text{Op\, Hom}}
(\Cal O_X, \Cal O_{X_{\text{red}}}).$$
In other words, we set $\Cal D_{\to}(n)$ to be the sheaf of
$n$--linear differential operators $\Cal O_X\times \dots \times \Cal O_X
\to \Cal O_{X_{\text{red}}}$ in the obvious sense. Clearly $\Cal D_{\to}$
is a $(\Cal D_{X_{\text{red}}}, \Cal D_X)$--bimodule. Similarly, let
$$\Cal D_{\leftarrow} = \{D_{\leftarrow}(n)\} = \underline{\text{ 
Op\,Hom}}(\Cal O_{X_{\text{red}}}, \Cal O_X)$$
be the collection of sheaves of $n$--linear differential operators
$\Cal O_{X_{\text{red}}}\times ... \times \Cal O_{X_{\text{red}}}
\to \Cal O_X$. This is a $(\Cal D_X, \Cal D_{X_{\text{red}}})$--bimodule.
The relative plethysm with these bimodules defines functors from 
$\Cal D_X$--SAlg to  $\Cal D_{X_{\roman{red}}}$--SAlg and back.
We claim that these functors are mutually inverse. More precisely,
we  have natural morphisms of sheaves of operads (in SVect)
$$\Cal D_{\leftarrow}\circ_{\Cal D_{X_{\text{red}}}}\Cal D_{\to}
\to \Cal D_{X_{\roman{red}}}, \quad \Cal D_{\to} \circ_{\Cal D_X} 
\Cal D_{\leftarrow} \to \Cal D_X,\leqno (2.3.3)$$
obtained by sheafification of (1.10.1--2).
In virtue of (1.5.7), it is enough to show that the morphisms
(2.3.3) are isomorphisms. To verify this, we can work locally,
over an affine open set $U\i X_{\text{red}}$. In this case
the supermanifold $(U, \Cal O_X|_U)$ is split, so 
Proposition 2.3.1(b) together with Theorem 1.10.3 imply
that our morphisms are isomorphisms after restruction on $U$.
This proves our theorem.

\newpage

\centerline{\bf References}

\medskip

[Ad] J.F. Adams. {\it Infinite Loop Spaces,} Princeton Univ. Press, 1978. 

\smallskip

[BJT] H.--J. Baues, M.~Jibladze, A.~Tonks. {\it Cohomology of monoids
in monoidal categories.} In: Proc. of Renaissance Conf., Contemp. Math.,
vol. 202, 137--165, Amer. Math. Soc., 1997.

\smallskip

[Bo] A. Borel et al. {\it Algebraic D-modules}, Academic Press, 1987. 

\smallskip 

[BG] A.~Beilinson, V.~Ginzburg. {\it Infinitesimal structure
of moduli spaces of $G$--bundles.} Int. Math. Res. Notes,
4 (1992), 63--74.

\smallskip

[CL] F. Chapoton, M. Livernet, {\it Pre-Lie algebras and the rooted
trees operad}, preprint math.QA/0002069. 

\smallskip 

[De] P.~Deligne. {\it Categories Tannakiennes}, in: Grothendieck Festschrift
(Eds. P. Cartier et al.) vol. II, p. 111--195, Birkhauser, Boston, 1990. 

\smallskip

[F] B.~Fresse, {\it Lie theory for formal groups over an operad}, 
J. of Algebra, 202 (1998), 455--511.

\smallskip

[Ger] M.~Gerstenhaber. {\it The cohomology structure of an associative
ring.} Ann. Math. 78 (1963), 59--103. 

\smallskip

[Ge--K] E.~Getzler, M.~Kapranov. {\it Cyclic operads and cyclic
homology.} In: ``Geometry, Topology, and Physics for Raoul,'' ed. by B. Mazur,
Cambridge, 1995, 167--201. 

\smallskip

[Gi--K] V.~Ginzburg, M.~Kapranov. {\it Koszul duality
for operads.} Duke Math.~J., 76 (1995), 203--272.

\smallskip

[J] A.~Joyal. {\it Foncteurs analytiques et esp\`eces de structures.}
Springer Lecture Notes in Math., vol. 1234, 126--159,
Springer Verlag 1986.

\smallskip

[Ka] M. Kashiwara. {\it Algebraic study of systems of partial differential
equations,} M\'emoires Soc. Math. France, 63 (1995), 1-72. 

\smallskip

[Ke] S.~Keel. {\it Intersection theory of moduli spaces of stable
$n$--pointed curves of genus zero.} Trans. AMS, 330 (1992), 545--574.

\smallskip

[Kn] F.~Knudsen. {\it Projectivity of the moduli space of stable
curves, II: the stacks $M_{g,n}.$} Math. Scand. 52 (1983),
161--199.
\smallskip

[Ku] B. A. Kupershmidt. {\it Non-Abelian phase spaces,}
J. Phys. A, 27 (1994), 2801-2810. 

\smallskip

[L] J.-L. Loday. {\it Op\'erations sur l'homologie cyclique des alg\'ebres
commutatives,} Invent. Math. 96 (1989), 205-230.

\smallskip

[Man1] Y.~I.~Manin. {\it Gauge Fields and Complex Geometry.}
Springer Verlag, 1988, 2nd edition 1997.

\smallskip

[Man2]  Y.~I.~Manin. {\it Frobenius manifolds, quantum cohomology, and moduli
spaces.} Colloq. Publ. Series, AMS, 1999.

\smallskip 

[Mar1] M.~Markl. {\it Models for operads.} Comm. in Algebra, 24 (1996),
1471--1500. 

\smallskip

[Mar2] M.~Markl. {\it Simplex, associahedron and cyclohedron,}
in:  Higher Homotopy Structures in Topology and Mathematical Physics,
J. McCleary Ed. (Contemporary Math. 227), p. 235--266.

\smallskip

[May] J.~P.~May. {\it The geometry of iterated loop spaces.}
Springer Lecture Notes in Math., vol. 271, Springer Verlag, 1972.

\smallskip

[P] I.~B.~Penkov. {\it D--modules on supermanifolds}, Inv. Math. 
71 (1983), 501--512.

\smallskip

[Pi] T. Pirashvili. {\it Hodge decomposition for higher order
Hochschild homology,} Ann. Sci. \'Ecole Norm. Sup. 33 (2000), 151-179. 

\smallskip

[Re] C. Rezk. {\it Spaces of algebra structures and cohomology of operads.}
Thesis, MIT, 1996.

\smallskip

[Ri] J.~F.~Ritt. {\it Differential algebra.} Amer. Math. Soc. 1950. 

\smallskip

[SKK] M. Sato,  M. Kashiwara, T. Kawai. {\it Hyperfunctions
and microdifferential equations}. Springer Lecture Notes in Math., vol.
287 (1973), p. 265--529. 

\smallskip

[Se] G.~B.~Segal. {\it Categories and homology theories.} Topology,
13 (1974), 293--312.

\smallskip

[Sm] V.~A.~Smirnov. {\it Homotopy theory of coalgebras.}
Math. USSR Izv., 27 (1986), 575--592.  

\smallskip

[W] N.~R.~Wallach. {\it Classical invariant theory and Virasoro algebra.}
In: ``Vertex Operators in Mathematics and Physics", J. Lepowsky,
S. Mandelstam, I.M. Singer, Eds. (MSRI Publications, Vol.3),  
475--482, Springer Verlag, 1979.

\enddocument